\documentclass{elsart}
\usepackage{amssymb,amsfonts}

\newcommand{\diag}{\mathop{\rm Diag}\nolimits}

\newcommand{\id}{\mathbb{I}}

\newcommand{\be}{\begin{equation}}
\newcommand{\ee}{\end{equation}}
\newcommand{\bea}{\begin{eqnarray}}
\newcommand{\eea}{\end{eqnarray}}
\newcommand{\beas}{\begin{eqnarray*}}
\newcommand{\eeas}{\end{eqnarray*}}

\newtheorem{theorem}{Theorem}
\newtheorem{lemma}{Lemma}

\newcount\minute
\newcount\hour
\def\currenttime{%
    \minute\time
    \hour\minute
    \divide\hour60
    \the\hour:\multiply\hour60\advance\minute-\hour\the\minute}
\begin{document}
\begin{frontmatter}
\title{A norm inequality for pairs of commuting positive semidefinite matrices}
\author{Koenraad M.R.\ Audenaert}
\address{
Department of Mathematics,
Royal Holloway University of London, \\
Egham TW20 0EX, United Kingdom \\[1mm]
Department of Physics and Astronomy, University of Ghent, \\
S9, Krijgslaan 281, B-9000 Ghent, Belgium}
\ead{koenraad.audenaert@rhul.ac.uk}
\date{\today, \currenttime}
\begin{abstract}
For $k=1,\ldots,K$, let $A_k$ and $B_k$ be positive semidefinite matrices such that, for each $k$, $A_k$ commutes with $B_k$.
We show that, for any unitarily invariant norm,
\[
|||\sum_{k=1}^K A_kB_k||| \le ||| (\sum_{k=1}^K A_k)\;(\sum_{k=1}^K B_k)|||.
\]
\end{abstract}

\end{frontmatter}
In this paper, we denote the vectors of eigenvalues and singular values of a 
matrix $A$ by $\lambda(A)$ and $\sigma(A)$, respectively. 
We adhere to the convention to sort singular values, and eigenvalues as well whenever they are real, in non-increasing order.
In general, for a real vector $x$, we will write $x^\downarrow$ for the vector with the same components as $x$ but sorted in 
non-increasing order.

For real $n$-dimensional vectors $x$ and $y$, we say that $x$ is \emph{weakly majorised} by $y$, denoted $x\prec_w y$, 
if and only if for $k=1,\ldots,n$,
$\sum_{i=1}^k x^\downarrow_i \le \sum_{i=1}^k y^\downarrow_i$. We say that $x$ is \emph{majorised} by $y$, denoted $x\prec y$, 
if and only if $x\prec_w y$ and $\sum_{i=1}^n x_i = \sum_{i=1}^n y_i$.
If, moreover, $x$ and $y$ are non-negative, we say that $x$ is weakly \emph{log-majorised} by $y$, 
denoted $x\prec_{w,\log} y$, if and only if for $k=1,\ldots,n$,
$\prod_{i=1}^k x^\downarrow_i \le \prod_{i=1}^k y^\downarrow_i$.

According to Weyl's Majorant Theorem (\cite{bhatia} Theorem II.3.6, or \cite{zhan}, Theorem 2.4), 
the vector of singular values of any matrix log-majorises
the vector of the absolute values of its eigenvalues: $|\lambda(A)|\prec_{\log}\sigma(A)$.
As $x\prec_{w,\log} y$ implies $x^r\prec_w y^r$ for any $r>0$, Weyl's Majorant Theorem can in slightly weaker form 
be stated as
\be
|\lambda(A)|^r\prec_{w}\sigma^r(A).\label{eq:WMT}
\ee

The sum of the $k$ largest singular values of a matrix defines a norm, known as the $k$-th Ky Fan norm.
The convexity of the Ky Fan norms can be expressed as a majorisation relation: for any $p$ such that $0\le p\le 1$,
\[
\sigma(pA+(1-p)B) \prec_w p\sigma(A)+(1-p)\sigma(B).
\]
When $A$ and $B$ are positive semidefinite, their singular values coincide with their eigenvalues and we have
\be
\lambda(pA+(1-p)B) \prec p\lambda(A)+(1-p)\lambda(B). \label{eq:convex}
\ee

For positive semidefinite matrices $A$ and $B$, the eigenvalues of $AB$ are real and non-negative. Furthermore
$\lambda(AB)\prec_{\log} \lambda(A)\lambda(B)$ (\cite{zhan} eq.~(2.5)). Hence, we also have
\be
\lambda(AB)\prec_{w} \lambda(A)\lambda(B).\label{eq:lamAB}
\ee

\bigskip

We start with a rather technical result concerning a majorisation relation for singular values:
\begin{lemma}\label{th:main}
Let $S$ be a general $n\times m$ complex matrix, and $L$ and $M$ two diagonal, positive semidefinite $m\times m$ matrices.
Then
\be
\sigma(SL\diag(S^*S)MS^*)\prec_w \sigma(SLS^*SMS^*).\label{eq:main}
\ee
\end{lemma}
\textit{Proof.}
We first show that
\be
\sigma(SL\diag(S^*S)MS^*)\prec_w \sigma(S(LM)^{1/2}S^*S(LM)^{1/2}S^*).\label{eq:1}
\ee
Since $L$, $M$, and $\diag(S^*S)$ are diagonal, they commute, and we can write
$SL\diag(S^*S)MS^* = S(LM)^{1/2}\diag(S^*S)(LM)^{1/2}S^*$. This is a positive semidefinite matrix, hence its singular
values are equal to its eigenvalues. The same is true for the right-hand side of (\ref{eq:1}).
Let us introduce $X=S(LM)^{1/4}$ and $T=X^*X\ge0$. 
Then we have to show $\lambda(X\diag(X^*X)X^*) \prec \lambda(XX^*XX^*)$, or
$\lambda(T\diag(T)) \prec \lambda(T^2)$.
Now note that there exist $J$ unitary matrices $U_j$ such that $\diag(T)=\sum_{j=1}^J (U_j T U_j^*)/J$.
Exploiting (\ref{eq:convex}) and inequality (\ref{eq:lamAB}) in turn, we obtain
\beas
\lambda(T\diag(T)) &=& \lambda(T^{1/2}\sum_{j=1}^J \frac{1}{J}(U_j T U_j^*) T^{1/2}) \\
&\prec& \sum_{j=1}^J \frac{1}{J} \lambda(T^{1/2} U_j T U_j^* T^{1/2}) \\
&=& \sum_{j=1}^J \frac{1}{J} \lambda(T U_j T U_j^* ) \\
&\prec_w& \sum_{j=1}^J \frac{1}{J} \lambda(T) \lambda(U_j T U_j^*) \\
&=& \sum_{j=1}^J \frac{1}{J} \lambda^2(T) = \lambda(T^2),
\eeas
which proves (\ref{eq:1}).

\medskip

Secondly, we show that
\be
\sigma(S(LM)^{1/2}S^*S(LM)^{1/2}S^*)\prec_w \sigma(SLS^*SMS^*).\label{eq:2}
\ee
Since $(LM)^{1/2}$ and $S^*S$ are both positive semidefinite, their matrix product has real, non-negative eigenvalues.
Thus,
\beas
\lambda^2((LM)^{1/2}S^*S)
&=& |\lambda(L^{1/2}S^*SM^{1/2})|^2 \\
&\prec_w& \sigma^2(L^{1/2}S^*SM^{1/2}),
\eeas
by Weyl's Majorant Theorem (eq. (\ref{eq:WMT}) with $r=2$).
This implies (\ref{eq:2}):
\beas
\sigma(S(LM)^{1/2}S^*S(LM)^{1/2}S^*) &=& \lambda((LM)^{1/2}S^*S(LM)^{1/2}S^*S) \\
&=& \lambda^2((LM)^{1/2}S^*S) \\
&\prec_w& \sigma^2(L^{1/2}S^*SM^{1/2}) \\
&=& \lambda^2((M^{1/2}S^*SLS^*SM^{1/2})^{1/2}) \\
&=& \lambda(M^{1/2}S^*SLS^*SM^{1/2}) \\
&=& \lambda(SLS^*SMS^*) \\
&=& |\lambda(SLS^*SMS^*)| \\
&\prec_w& \sigma(SLS^*SMS^*),
\eeas
where in the last line we again exploited Weyl's Majorant Theorem (eq. (\ref{eq:WMT}) with $r=1$).

Combining (\ref{eq:1}) with (\ref{eq:2}) yields (\ref{eq:main}).
\qed

\bigskip

We can now state and prove the main result of this paper.
\begin{theorem}\label{cor:1}
For $k=1,\ldots,K$, let $A_k$ and $B_k$ be positive semidefinite $d\times d$
matrices such that, for each $k$, $A_k$ commutes with $B_k$.
Then for all unitarily invariant norms
\be
|||\sum_{k=1}^K A_kB_k||| \le ||| (\sum_{k=1}^K A_k)\;(\sum_{k=1}^K B_k)|||.\label{eq:cor1}
\ee
\end{theorem}
\textit{Proof.}
Let $A_k$ and $B_k$ have eigenvalue decompositions 
\[
A_k = U_k a_k U_k^*,\quad  B_k =U_k b_k U_k^*,
\]
where the $U_k$ are unitary matrices, and $a_k$ and $b_k$ are positive semidefinite diagonal matrices.
Let 
\[
L=\bigoplus_{k=1}^K a_k, \quad M=\bigoplus_{k=1}^K b_k,\quad S=(U_1 | U_2 | \cdots | U_K).
\]
Then 
\[
\sum_{k=1}^K A_k = SLS^*,\quad \sum_{k=1}^K B_k = SMS^*, \quad \sum_{k=1}^K A_kB_k = SLMS^*.
\]
In addition, $\diag(S^*S)=\id$ since all columns of $S$ are normalised.
By Lemma \ref{th:main}, we then have
\[
\sigma(\sum_{k=1}^K A_kB_k)\prec_w \sigma((\sum_{k=1}^K A_k)\;(\sum_{k=1}^K B_k))
\]
which is equivalent to (\ref{eq:cor1}).
\qed

A simple consequence of Theorem \ref{cor:1} is that for any set of $K$ positive semidefinite matrices $A_k$,
all positive functions $f$ and $g$, and all unitarily invariant norms,
\be
|||\sum_{k=1}^K f(A_k)g(A_k)||| \le ||| (\sum_{k=1}^K f(A_k))\;(\sum_{k=1}^K g(A_k))|||.
\ee
Setting $K=2$, $f(x)=x^p$ and $g(x)=x^q$ yields the inequality
\be
|||A^{p+q} +B^{p+q}||| \le ||| (A^p+B^p)(A^q+B^q)|||,
\ee
which was recently conjectured by Bourin \cite{bourin}.

\section*{Acknowledgments}
We acknowledge support by an Odysseus grant from the Flemish FWO.



\begin{thebibliography}{99}
\bibitem{bhatia} R.~Bhatia, \textit{Matrix Analysis}, Springer, Berlin (1997).
\bibitem{bourin} J.-C.~Bourin, ``Matrix subadditivity inequalities and block-matrices'', 
Internat.\ J.\ Math.\ \textbf{20}, 679--691 (2009).
\bibitem{zhan} X.~Zhan, \textit{Matrix Inequalities}, Springer, Berlin (2002).
\end{thebibliography}
\end{document}